\documentclass[a4paper,article,11pt,fleqn,oneside]{memoir} 
\setsecnumdepth{section}
\usepackage {amsmath,amsthm,amsxtra,amssymb}
\usepackage[fleqn,tbtags]{mathtools}      
\usepackage{upgreek} 

\usepackage[utf8]{inputenc}
\usepackage[english]{babel}
\usepackage{lmodern}   
\usepackage[T1]{fontenc}	 
\usepackage {bm} 
\usepackage{fix-cm} 
\usepackage{microtype}	
\usepackage{float} 
\usepackage[compatibility=false]{caption}

\checkandfixthelayout[nearest]

\usepackage[babel=true]{csquotes}
\usepackage{graphicx}
\usepackage[ansinew]{} 
\usepackage{fixltx2e,mparhack}
\usepackage[fixpdftex]{xcolor} 
\usepackage{tikz}
\usepackage{pgf}
\usetikzlibrary{fit} 
\usetikzlibrary[shapes.geometric]
\usetikzlibrary{matrix,arrows,petri,topaths,positioning,automata}

\setlrmarginsandblock{3.14cm}{3.5cm}{*}
\setulmarginsandblock{3.14cm}{3.3cm}{*}
\checkandfixthelayout

\headstyles{memman}
\chapterstyle{article}
\setsubsecheadstyle{\centering\normalfont}
\setbeforesubsecskip{-1\onelineskip plus -1ex minus -.2ex}
\setaftersubsecskip{1\onelineskip plus .2ex}
\usepackage[colorlinks=true, pdfstartview=FitV, linkcolor=blue, citecolor=blue, urlcolor=blue]{hyperref}

\let \rho   = \varrho
\let \phi   = \varphi

\let \Phi = \varPhi
\let \Psi = \varPsi
\let \Gamma = \varGamma
\let \Theta = \varTheta

\DeclareMathOperator{\Set}{\textbf{Set}}

\DeclareMathOperator{\id}{id}
\DeclareMathOperator{\FGR}{\textbf{Graph}_{\textbf{\F}}}

\DeclareMathOperator{\pb}{pb}

\DeclareMathOperator{\col}{Col}

\DeclareMathOperator{\Cc}{C}
\DeclareMathOperator{\C}{C}
\DeclareMathOperator{\Hom}{H}

\DeclareMathOperator{\Sim}{S}
\DeclareMathOperator{\Graph}{Graph}
\DeclareMathOperator{\Pat}{Pat}

\newcommand {\F}{F} 
\newcommand {\Fb}{\bar{F}} 
\newcommand{\Uu}{U}

\newcommand {\G}{G} 
\newcommand {\Kk}{\mathcal{K}}

\newcommand {\Pp}{\mathcal{P}}

\newcommand{\ug}[1][]{#1_{\scriptscriptstyle{g}}}
\newcommand{\vv}[1][]{#1_{\scriptscriptstyle{v}}}
\newcommand{\ee}[1][]{#1_{\scriptscriptstyle{e}}}

\newcommand{\ui}[1][]{#1_{i\in I}}

\newcommand{\inkl}[2]{\iota_{(\scriptstyle{#1}{,}{\scriptstyle{#2}})}}

\newcommand {\he}{^{\scriptscriptstyle{(1)}}}   
\newcommand {\hz}{^{\scriptscriptstyle{(2)}}}

\newcommand {\hi}{^{\scriptscriptstyle{(i)}}}

\definecolor{rot}{rgb}{1,0,0}
\definecolor{weis}{rgb}{1,1,1}

\theoremstyle{definition}
\newtheorem{defi}{Definition}[section]
\newtheorem{numbsp}[defi]{Example}
\theoremstyle{remark}
\newtheorem{numbem}[defi]{Remark}

\theoremstyle{plain}

\newtheorem{satz}[defi]{Theorem}
\newtheorem{prop}[defi]{Proposition}
\newtheorem{theorem}[defi]{Theorem}

\begin{document}

\title{A coalgebraic model of graphs}
\author{Christian J\"akel, \textit{christian.jaekel(at)tu-dresden.de}\\
Technische Universit\"at Dresden}
\date{January 2016}
\maketitle
\begin{abstract}
In this note, we model various types of graphs, relational systems and multisets as coalgebras over $\Set\times \Set$ and use the theory of coalgebras over arbitrary categories to conclude properties of the category of graphs. This point of view forces the formulation of a Co-Birkhoff like theorem for graphs.
\end{abstract}

\section{Introduction}
For a covariant endofunctor $\F$ in $\Set$, we can model various types of graphs through comma categories over $\Set$. This yields a triple $\G=(V,E,g)$, where $V$ is the \emph{vertex set}, $E$ the \emph{edge set} and $g:E\rightarrow\F V$ the \emph{strucure map} of the $\F$-graph $\G$. For two $\F$-graphs $\G\he=(V\he,E\he,g\he)$ and $\G\hz=(V\hz,E\hz,g\hz)$, a homomorphism is a pair of maps $\phi=(\phi\vv,\phi\ee)$, where $\phi\vv:V\he\rightarrow V\hz$ defines the \emph{vertex map} and $\phi\ee: E\he\rightarrow E\hz$, the \emph{edge map}, such that $g\hz\circ\phi\ee=\F(\phi\vv)\circ g\he$ holds.\\

Next, we list examples for choices of $F$ and resultant graph structures. 

\begin{numbsp}
\begin{itemize}
\item The identity functor yields a model for multisets, \textit{i.e.}, a triple $(\mathbb{N},E,ar)$, where the structure map $ar$ assigns an arity to each element of $E$.
\item By considering the functor $\F:V\mapsto \mathfrak{P}_{_{1,2}}V$, which assigns to a set its singleton and two-element subsets, we can model undirected graphs $g:E\rightarrow\mathfrak{P}_{_{1,2}}V$ and the induced homomorphisms are incidence preserving maps.
\item Directed graphs $g:E\rightarrow V\times V$ can be represented with  $\F:V\mapsto V\times V$ and homomorphisms preserve soure- and target nodes.
\item Let $\mathfrak{P}$ be the powerset functor. In this case, $g:E\rightarrow\mathfrak{P}V$ or $g:E\rightarrow V\times \mathfrak{P}V$ yields a model of hypergraphs or directed hypergraphs respectively. 
\item Hybrid graphs can be described by taking the sum of different type functors, as for instance, $\F V=(V\times V) + \mathfrak{P}_{_{1,2}} V $ for a model of graphs with directed and undirected edges. In this framework, relational systems fit in too \cite{FoniokRelsystems}.
\item Let $X\vv$ represent a set of vertex colors and $X\ee$ a set of edge colors. For every type functor $\F$, we define colored graphs $g:E\rightarrow X\ee\times\F(X\vv\times V)$ by the functor $\tilde{\F}(V):=X\ee\times\F(X\vv\times V)$. For example, with $X\vv=X\ee=[0,1]$, we get fuzzy graphs.
\end{itemize}
\end{numbsp}
Hence, various mathematical structures are generalized and can be described through $\F$-graphs.\\
In \cite{Ruth}, it was noticed that directed graphs can be considered as coalgebras over $\Set\times\Set$, by the functor $\Fb:(V,E)\mapsto(1,V\times V)$, together with the structure map $g:(V,E)\rightarrow (1,V\times V)$. The induced homomorphisms are source- and target node preserving.\\
Inspired by this example, we developed a purely coalgebraic model of various types of graphs. We will present it in this note and apply the theory of coalgebras over arbitrary base categories (as in \cite{HughesDiss,AdamekPorst}), to conclude structural results about the category of graphs. This yields a very economical approach to the category of graphs. Thereby, we generalize results of \cite{DeCatGraph,Williams}, as for example the construction of cofree graphs and results about the category of multisets \cite{IsahMultiset,SinghMultiset}. 

\section{Graphs As Coalgebras}
We will define graphs as coalgebras and conclude immediate consequences for the category of graphs.
\begin{defi} Let $\F:\Set\rightarrow\Set$ be a covariant endofunctor. We extend $\F$ to be a covariant endofunctor $\Fb:\Set\times\Set\rightarrow\Set\times\Set$, via $(V,E)\mapsto(1,\F V)$ and $(\phi\vv,\phi\ee)\mapsto(!,\F(\phi\vv))$ , where $1$ is the terminal object in $\Set$ and $!:V\rightarrow 1$ is the induced unique map.\\
An \emph{$\F$-graph} is a coalgebra $((V,E),(!,g))$, where $g: E\rightarrow \F V$ determines the graphs internal structure. A homomorphism $\phi=(\phi\vv,\phi\ee)$ is defined through the following commutative diagram.
\begin{figure}[H]
\begin{minipage}[h]{5cm}
\begin{center}
\begin{tikzpicture}[description/.style={fill=white,inner sep=2pt},>=stealth']
\matrix (m) [matrix of math nodes, row sep=2.71em,
column sep=3.14em, text height=1.5ex, text depth=0.25ex]
{ (V\he,E\he) & (V\hz,E\hz) \\
 (1,\F V\he) & (1,\F V\hz). \\ };
\path[->,font=\scriptsize](m-1-1) edge node[auto] {$ (\phi\vv,\phi\ee) $} (m-1-2);
\path[->,font=\scriptsize](m-1-1) edge node[auto,swap] {$ (!, g\he) $} (m-2-1);
\path[->,font=\scriptsize](m-2-1) edge node[auto,swap] {$ (!,\F(\phi\vv)) $} (m-2-2);
\path[->,font=\scriptsize](m-1-2) edge node[auto] {$ (!,g\hz)$} (m-2-2);
\end{tikzpicture}
\end{center}
\end{minipage}
\hspace{1cm}
\begin{minipage}[h]{1cm}
$\cong$
\end{minipage}
\hspace{0.15cm}
\begin{minipage}[h]{1cm}
\begin{center}
\begin{tikzpicture}[description/.style={fill=white,inner sep=2pt},>=stealth']
\matrix (m) [matrix of math nodes, row sep=2.71em,
column sep=3.14em, text height=1.5ex, text depth=0.25ex]
{ (V\he,E\he) & (V\hz,E\hz) \\
 \Fb(V\he,E\he) & \Fb(V\hz,E\hz) . \\ };
\path[->,font=\scriptsize](m-1-1) edge node[auto] {$ (\phi\vv,\phi\ee) $} (m-1-2);
\path[->,font=\scriptsize](m-1-1) edge node[auto,swap] {$  g\he $} (m-2-1);
\path[->,font=\scriptsize](m-2-1) edge node[auto,swap] {$ \Fb(\phi\vv,\phi\ee) $} (m-2-2);
\path[->,font=\scriptsize](m-1-2) edge node[auto] {$ g\hz$} (m-2-2);
\end{tikzpicture}
\end{center}
\end{minipage}
\end{figure}
For a simplified notation, we will write $(V,E,g)$, instead of $((V,E),(!,g))$. The category of graphs will be denoted by $\FGR$. It is equipped with the forgetful functor $\Uu:\FGR\rightarrow\Set\times\Set$.
\end{defi}
By applying the definition of $\Fb$, the following can be shown straight forward.
\begin{prop} The extension $\Fb$ preserves all limits that $\F$ preserves and $\Fb$ is bounded iff $\F$ is bounded. 
\end{prop} 
Next, we will list structural results for $\FGR$, which follow immediately from the theory of coalgebras. Thereby, we also use the fact that $\Set\times\Set$ inherits most of its properties from $\Set$.

\begin{prop}
The forgetful functor $\Uu$ creates colimits. Hence, colimits in $\FGR$ are constructed as pairs of colimits in $\Set$. Furthermore, $\Uu$ creates any limit preserved by $\F$ (see \cite[Theorem 1.2.4, Theorem 1.2.7]{HughesDiss}).
\end{prop} 

\begin{prop}
Isomorphisms are bijective\footnote{By injective, surjective and bijective maps $\phi=(\phi\vv,\phi\ee)$, we mean that $\phi\vv$ and $\phi\ee$ are injective, surjective and bijective respectively.} maps $\phi=(\phi\vv,\phi\ee)$ and epimorphisms are the surjective ones (see \cite[4.18 Corollary]{AdamekPorst}). Regular monomorphisms are the injective maps (see \cite[Theorem 3.4]{GummWeakLimit}). 
\end{prop}

\begin{prop}
Subgraphs are defined as regular subobjects (see \cite[Definition 2.2.1]{HughesDiss}) and the category $\FGR$ has an epi-regular mono factorization system, which is created by $\Uu$ (see \cite[4.23 Remark]{AdamekPorst}).
\end{prop}

\begin{prop}\label{EqualizerProp} The category $\FGR$ has all equalizers and the equalizer object is cogenerated from the respective equalizer object in $\Set\times\Set$ (see \cite[Theorem 2.4.1]{HughesDiss}).
\end{prop}
The cofree graph for a color set $X$ can be constructed as follows: 
\begin{satz}
For a set of colors $X=(X\vv,X\ee)$, the cofree graph over $X$ is given through $(X\vv,X\ee\times\F X\vv,\pi_{\F X\vv})$ together with the map $\varepsilon_X=(\id_{X\vv},\pi_{X\ee})$, where $\pi_{\F X\vv}$ and $\pi_{X\ee}$ are the canonical projections of the product $X\ee\times\F X\vv$.
\begin{proof}
We construct the transfinite cochain as in \cite[2.23 Cofree Coalgebra Construction]{AdamekPorst} and see that it stops after two steps. We define: $X^0_\#=(1,1)$. Then,
\begin{align*}
{X^1_\#}&=(X\vv,X\ee)\times\Fb(1,1)=( X\vv\times 1, X\ee\times\F 1)=(X\vv, X\ee\times\F 1),\\
X^2_\#&=(X\vv,X\ee)\times\Fb(X\vv,X\ee\times\F 1)=(X\vv,X\ee)\times(1,\F X\vv)=(X\vv,X\ee\times\F X\vv),\\
X^3_\#&=(X\vv,X\ee)\times\Fb(X\vv,X\ee\times\F X\vv)=(X\vv,X\ee)\times(1,\F X\vv)=(X\vv,X\ee\times\F X\vv).
\end{align*}
\end{proof}
\end{satz}

The above construction gives rise to a right adjoint functor $\Uu\dashv\Cc$ with $\Cc:\Set\times\Set\rightarrow\FGR$. Consequently, it holds that for every graph $\G$ and every coloring $\gamma:\Uu\G\rightarrow(X\vv,X\ee)$ there exists a unique homomorphism $\overline{\gamma}$, such that $\gamma=\varepsilon_X\circ\Uu(\overline{\gamma})$. 
\begin{numbsp}
We consider undirected graphs of type $\mathfrak{P}_{_{1,2}}$. Let $X\vv=\{r,g\}$ and $X\ee=\{1,2\}$. The resulting cofree graph is pictured below (left-hand), together with the graph coloring $\gamma$ and the assigned colors are labels on the graph $G$ (right-hand). 
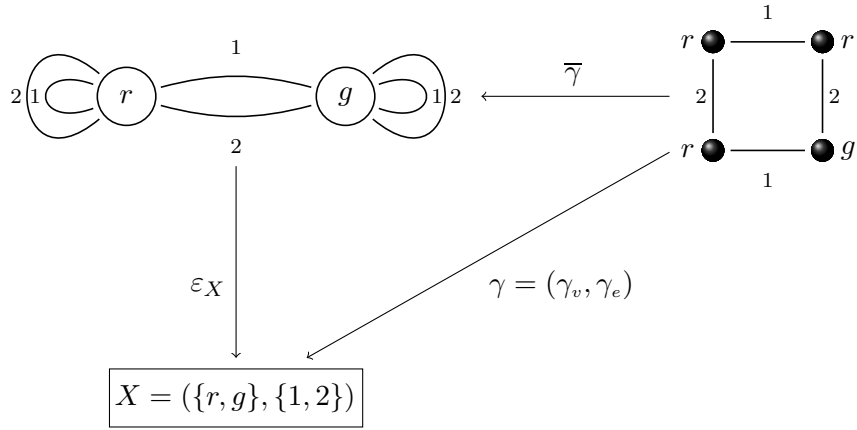
\begin{figure}[H]
\begin{center}
\begin{tikzpicture}[latticeelem/.style={rectangle,text centered,draw=white,inner sep=2pt,minimum height=2em},latticeelemblack/.style={rectangle,text centered,draw=black,inner sep=2pt,minimum height=2em},thin]

 \node[latticeelem,anchor=center](CofreeGraph)
 {
\begin{tikzpicture}[scale=0.8,descr/.style={fill=white,inner sep=2.5pt},>=stealth',shorten >=2pt , shorten <=2pt,node distance=2.8cm,semithick,scale=0.6,every fit/.style={ellipse,draw,inner sep=-2pt},every loop/.style={},cross line/.style={preaction={draw=white,-,line width=6pt}}]
   	
 	  \node[draw,circle] (a) at (0,0){$r$};  
 	  \node[draw,circle] (b) at (6,0){$g$}; 	   	
 	 
	\path[-,font=\scriptsize](a) edge[bend left=15] node[auto]{$1$}(b);
 	\path[-,font=\scriptsize](a) edge[bend right=15] node[auto,swap]{$2$}(b); 	 	
	\path[-,font=\scriptsize](a) edge[loop left,in=200,out=160,min distance=20mm] node[auto,swap]{$1$}(a);
	\path[-,font=\scriptsize](a) edge[loop left,in=220,out=140,min distance=36mm] node[auto,swap]{$2$}(a);	
	\path[-,font=\scriptsize](b) edge[in=340,out=20,min distance=20mm] node[auto]{$1$} (b);
	\path[-,font=\scriptsize](b) edge[in=320,out=40,min distance=36mm] node[auto]{$2$} (b);
 	 
\end{tikzpicture}
  };  
	
\node[latticeelem,right of=CofreeGraph,node distance=7cm,anchor=center] (TestGraph)
 {

\begin{tikzpicture}[scale=0.8,>=stealth',shorten >=2pt , shorten <=2pt,auto,node distance=2.8cm,
                    semithick,scale=0.6,every fit/.style={ellipse,draw,inner sep=-2pt}]
   	
 	  \node[minimum size=0.3cm,circle,ball color=black,label=left:{$r$}](a) at (0,0){};  
 	  \node[minimum size=0.3cm,circle,ball color=black,label=right:{$r$}] (b) at (3,0){};
 	  \node[minimum size=0.3cm,circle,ball color=black,label=right:{$g$}] (c) at (3,-3){};
		\node[minimum size=0.3cm,circle,ball color=black,label=left:{$r$}] (d) at (0,-3){}; 
 	 
 	 \path[-,font=\scriptsize](a) edge node[auto]{$1$}(b);
 	 \path[-,font=\scriptsize](b) edge node[auto]{$2$}(c);
 	 \path[-,font=\scriptsize](c) edge node[auto]{$1$}(d);
 	 \path[-,font=\scriptsize](d) edge node[auto]{$2$}(a);
	\path[-,font=\scriptsize](a) edge[white, loop above] (a);
\end{tikzpicture}  
  };
\node[latticeelemblack,below of=CofreeGraph,node distance=4cm,anchor=center] (ColorSet)
 {
   $X=(\{r,g\},\{1,2\})$
   
  };
\path[->,shorten >=8pt](TestGraph)   edge  node[auto] {$\gamma=(\gamma\vv,\gamma\ee)$}(ColorSet);
\path[->,shorten >=4pt,shorten <=-15pt](CofreeGraph) edge  node[auto,swap] {$\varepsilon_X$} (ColorSet);
\path[->,shorten >=2pt](TestGraph)  edge  node[auto,swap] {$\overline{\gamma}$}(CofreeGraph);  
\end{tikzpicture}
\end{center}
\end{figure}
\end{numbsp}

\begin{numbem}
For a graph $\G=(V,E,g)$ and a coloring $\gamma:\Uu\G\rightarrow(X\vv,X\ee)$, we can define colored graphs with respect to the functor 
$(X\vv,X\ee)\times\Fb(V,E)\cong(X\vv,X\ee\times\F V)$. The structure map of the colored graph is given as $(\gamma\vv,(\gamma\ee,g)):(V,E)\rightarrow(X\vv,X\ee\times\F V)$.
\end{numbem}

\begin{numbem}
For $\F=\mathfrak{P}$, the extension $\bar{\mathfrak{P}}$ yields a non-accessible covarietor.
\end{numbem}

Using the dual of \cite[4.13 Theorem]{AdamekPorst}, we can conclude that $\FGR$ is complete. Additionally to this fact, we will present an explicit construction of products, which does, together with proposition \ref{EqualizerProp}, also imply the completeness of $\FGR$.

\begin{theorem}\label{products} 
Let $(\G\hi=(V\hi,E\hi,g\hi))\ui$ be a family of $\F$-graphs. Their product is given as $\prod\G\hi=(\prod V\hi,E_{prod},g_{prod})$, where $E_{prod}$ and $g_{prod}:=\pb(\alpha)$ are defined through the following pullback square. 
 \begin{center}
\begin{tikzpicture}[state/.style={rectangle,rounded corners,draw=white,minimum height=2em,
           inner sep=2pt,text centered},->,>=stealth',descr/.style={fill=white,inner sep=2.5pt},scale=0.9]  
\node[state] (a){$E_{prod}$};
\node[state, right of=a,node distance=3cm,anchor=center](c){$\F(\prod V\hi)$};
\node[state, below of=a,node distance=2.1cm,anchor=center](d){$\prod E\hi$};
\node[state, right of=d,node distance=3cm,anchor=center](e){$\prod\F V\hi$};
\path[->,font=\scriptsize](a) edge node[auto] { $\pb(\alpha)$ } (c);
\path[->,font=\scriptsize](a) edge node[auto,swap] { $\pb(\beta)$ } (d);
\path[->,font=\scriptsize](c) edge node[auto] { $\beta$ } (e);
\path[->,font=\scriptsize](d) edge node[auto,swap] {$\alpha$ } (e);
\end{tikzpicture}
\end{center}
The projection homomorphisms are $\pi\hi:=(\pi\hi\vv,\pi\hi\ee\circ\pb(\beta)):\prod\G\hi\rightarrow\G\hi$, where $\pi\hi\vv$ and $\pi\hi\ee$ are the projections of $\prod V\hi$ and $\prod E\hi$ respectively. 
\begin{proof}
It is straight forward to show that the $\pi\hi$ are homomorphisms. Next, we notice that every $\F$-graph $\G\hi=(V\hi,E\hi,g\hi)$ is a subgraph of $\Cc(V\hi,E\hi)$ (see \cite[Theorem 2.1.16]{HughesDiss}) and $\prod\Cc(V\hi,E\hi)=\Cc(\prod V\hi,\prod E\hi)=(\prod V\hi,\prod E\hi\times\F(\prod V\hi))$ holds. According to \cite[Theorem 2.4.3]{HughesDiss}, the product of the $\G\hi$ is cogenerated, with respect to $\prod\Cc(V\hi,E\hi)$, by the following pullback in $\Set\times\Set$.
\begin{center}
\begin{tikzpicture}[state/.style={rectangle,rounded corners,draw=white,minimum height=2em,
           inner sep=2pt,text centered},->,>=stealth',descr/.style={fill=white,inner sep=2.5pt},transform shape,scale=1]  
\node[state] (a){$(\prod V\hi,\tilde{E})$};
\node[state, right of=a,node distance=4.8cm,anchor=center](b){$(\prod V\hi,\prod E\hi\times\F(\prod V\hi))$};
\node[state, below of=a,node distance=2.1cm,anchor=center](c){$(\prod V\hi,\prod E\hi)$};
\node[state, right of=c,node distance=4.8cm,anchor=center](d){$(\prod V\hi,\prod E\hi\times\prod\F V\hi)$};
\path[->,font=\scriptsize](a) edge node[auto,swap] { $$ } (b);
\path[->,font=\scriptsize](b) edge node[auto,swap] { $$ } (d);
\path[->,font=\scriptsize](a) edge node[auto,swap] { $$ } (c);
\path[->,font=\scriptsize](c) edge node[auto,swap] { $$ } (d);
\end{tikzpicture}
\end{center}
By diagram chasing, one can show that $\tilde{E}\cong E_{prod}$. Hence, the largest subgraph contained in $(\prod V\hi,\tilde{E})$ is 
$(\prod V\hi,E_{prod},\pb(\alpha))$.
\end{proof}
\end{theorem} 
It follows that the edge set of the product is a subset of $E\he\times E\hz\times\F(V\he\times V\hz)$. This yields the usual product of undirected graphs.
\begin{numbsp}
In the product of two $\mathfrak{P}_{_{1,2}}$-graphs $\G$ and $\tilde{G}$, for each $e\in E$ with $g(e)=\{v,w\}$ and $\tilde{e}\in\tilde{E}$ with $\tilde{g}(e)=\{\tilde{v},\tilde{w}\}$, there are two edges in $\G\times\tilde{G}$, namely $(e,\tilde{e},\{(v,\tilde{v}),(w,\tilde{w})\})$ and $(e,\tilde{e},\{(v,\tilde{w}),(w,\tilde{v})\})$.
\end{numbsp}

At last, we define covarieties and coequation satisfaction in the usual way.
\begin{defi} 
Let $\G=(V,E,g)$ be an $\F$-graph and $\Cc X$ the cofree graph over a color set $X=(X\vv,X\ee)$. We define $\col_X(\G)=\{\gamma\mid\gamma:\Uu\G\rightarrow X\}$ to be the collection of all colorings of $\G$. A subset $P\subseteq(X\vv,X\ee)$ is called \emph{pattern} over $\Cc X$. We say that a pattern $P$ holds in $\G$ if for all $\gamma\in\col_X(\G)$, we have that the image $\overline{\gamma}[\G]$ is a subgraph of $[P]\ug$, the graph cogenerated by $P$. For that, we  write $\G\vDash P$. This means that for all colorings, the induced homomorphism $\overline{\gamma}$ factors through the largest in $P$ contained subgraph, or in other words that $\G$ is regular-projective with respect to the inclusion morphism $\inkl{[P]\ug}{\C X}$.
\begin{figure}[H]
\begin{center}
\begin{tikzpicture}[description/.style={fill=white,inner sep=2pt},>=stealth']
\matrix (m) [matrix of math nodes, row sep=3.5em,
column sep=2.5em, text height=1.3ex, text depth=0.25ex] 
{\G & & \Cc X\\
    & & {[}P{]}\ug \\ };
\path[->,font=\scriptsize](m-1-1) edge node[auto] {$ \overline{\gamma} $} (m-1-3);
\path[dashed,->,font=\scriptsize](m-1-1) edge node[auto,swap] {$ $} (m-2-3);
\path[right hook->,shorten <=4,font=\scriptsize](m-2-3) edge node[auto,swap] {$\inkl{[P]\ug}{\Cc X}$} (m-1-3);
\end{tikzpicture}
\end{center}
\end{figure}
Let $\Pp$ be a collection of patterns and $\Kk$ a class of $\F$-graphs. We say $\G\vDash \Pp$ if every $P\in \Pp$ holds in $\G$. Analogously, $\Kk\vDash P$ provided $\G\vDash P$ for every $\G\in\Kk$.
\end{defi}

\begin{satz} A class $\Kk$ of $\F$-graphs is a co-variety (\textit{i.e.}, closed under subgraphs, homomorphic images and coproducts) if and only if there is a collection $\Pp$ of patterns, such that for all $\G$, it holds that $\G\in\Kk\Leftrightarrow\G\vDash \Pp$ (see \cite[Theorem 3.6.21]{HughesDiss}).
\end{satz}

We define: $\Graph(\Pp):=\{\G\in\FGR\mid \G\vDash\Pp\}$ and $\Pat(\Kk):=\{P \mid \Kk\vDash P\}.$

\begin{satz}[Co-Birkhoff Theorem For $\F$-Graphs] For a class $\Kk$ of $\F$-graphs, we have: $\Sim\Hom\Sigma(\Kk)=\Graph(\Pat(\Kk)).$\footnote{Here $\Sim$ is the class of all subgraphs, $\Hom$ the class of all homomorphic images and $\Sigma$ the class of all coproducts of elemts from $\Kk$.}
\end{satz}

Furthermore, it holds that: 

\begin{satz}
If $\F$ is bounded, then every covariety can be defined through a set of coequations (see \cite[Theorem 3.7.22]{HughesDiss}).
\end{satz}

\section{Conclusion}
As for example groups or ring are algebras, we showed that various types of graphs and mathematical structures like multisets or relational systems, together with their structure preserving homomorphisms, can be considered as coalgebras. This yields new examples for coalgebras. Furthermore, from this point of view, the categorical properties of graphs follow immediately from the theory of universal coalgebras and do not have to be shown in a tedious treatment of special cases.


\begingroup
\renewcommand{\bibname}{References}
\setlength\bibitemsep{0pt}

\endgroup


\begin{thebibliography}{10}

\bibitem{AdamekPorst}
J.~Ad\'{a}mek and H.-E.~Porst, {\em On Varieties and Covarieties in a Category}.
\newblock Cambridge University Press, Mathematical. Structures in Comp. Sci., pp.~201--232, 2003.

\bibitem{FoniokRelsystems}
J.~ Foniok, {\em Homomorphisms and Structural Properties of Relational Systems}.
\newblock  PhD. thesis, Charles University, Prague, 2007.

\bibitem{GummWeakLimit}
H.~P. Gumm and T.~Schr\"{o}der, {\em Coalgebraic structure from weak limit
  preserving functors}. 
\newblock Coalgebraic Methods in Computer Science, 2000.

\bibitem{HughesDiss}
J.~Hughes, {\em A Study of Categories of Algebras and Coalgebras}.
\newblock PhD thesis, Carnegie Mellon University, 2001.

\bibitem{IsahMultiset}
A.~I. Isah and Y.~Tella, {\em The Concept of Multiset Category}.
\newblock British Journal of Mathematics {\&} Computer Science, 2015.

\bibitem{DeCatGraph}
D.~J. Plessas, {\em The Categories of Graphs}.
\newblock PhD thesis, The University of Montana Missoula, MT, 2011.

\bibitem{Ruth}
J.~Rutten, {\em Universal coalgebra: a theory of systems}.
\newblock Theoretical Computer Science,  pp.~3--80, 2000.

\bibitem{SinghMultiset}
D.~Singh, {\em A Note on Category of Multisets (MUL)}.
\newblock International Journal of Algebra, pp.~73--78, 2013.
  
\bibitem{Williams}
K.~K. Williams, {\em The category of graphs}.
\newblock master thesis, Texas Tech University, 1971.

\end{thebibliography}
\end{document}